\documentclass[12pt]{article}
\usepackage{mathrsfs}
\usepackage{amsfonts}
\usepackage{graphicx}
\usepackage{amsfonts, amsmath, amssymb, latexsym, epsfig}
\newcommand{\blst}{\begin{trivlist}}
\newcommand{\elst}{\end{trivlist}}

\newtheorem{thm}{Theorem}[section]
\newtheorem{prop}{Proposition}[section]
\newtheorem{cor}{Corollary}[section]
\newtheorem{lem}{Lemma}[section]
\newtheorem{conj}{Conjecture}[section]
\newtheorem{exa}{Example}[section]
\newtheorem{defn}{Definition}[section]
\newtheorem{rem}{Remark}[section]

\newcommand{\ben}{\begin{enumerate}}
\newcommand{\een}{\end{enumerate}}
\newcommand{\ble}{\begin{lem}}
\newcommand{\ele}{\end{lem}}
\newcommand{\bth}{\begin{thm}}
\renewcommand{\eth}{\end{thm}}
\newcommand{\bpr}{\begin{prop}}
\newcommand{\epr}{\end{prop}}
\newcommand{\bco}{\begin{cor}}
\newcommand{\eco}{\end{cor}}
\newcommand{\bcon}{\begin{conj}}
\newcommand{\econ}{\end{conj}}
\newcommand{\bde}{\begin{defn}}
\newcommand{\ede}{\end{defn}}
\newcommand{\bex}{\begin{exa}}
\newcommand{\eex}{\end{exa}}
\newcommand{\barr}{\begin{array}}
\newcommand{\earr}{\end{array}}
\newcommand{\btab}{\begin{tabular}}
\newcommand{\etab}{\end{tabular}}
\newcommand{\beq}{\begin{equation}}
\newcommand{\eeq}{\end{equation}}

\newcommand{\bea}{\begin{eqnarray*}}
\newcommand{\eea}{\end{eqnarray*}}
\newcommand{\beaa}{\begin{eqnarray}}
\newcommand{\eeaa}{\end{eqnarray}}
\newcommand{\bce}{\begin{center}}
\newcommand{\ece}{\end{center}}
\newcommand{\bpi}{\begin{picture}}
\newcommand{\epi}{\end{picture}}
\newcommand{\bfi}{\begin{figure} \begin{center}}
\newcommand{\efi}{\end{center} \end{figure}}

\newcommand{\bsl}{\begin{slide}{}}
\newcommand{\esl}{\end{slide}}

{}
{}
{}
{}
{}
{}
{}
{}
{}
{}
{}
{}
{}
{}
{}
{}
{}
{}
{}
{}
{}
{}
{}
{}
{}
{}
{}
{}
\newenvironment{proof}{
\par
\noindent {\bf Proof.}\rm}{\mbox{}\hfill\rule{0.5em}{0.809em}\par}
\begin{document}
\title{Circular Peaks and Hilbert Series}
\author{Pierre Bouchard$^{a}$
 \and Jun Ma$^{b}$\thanks{Email address of the corresponding author:majun@math.sinica.edu.tw}\\
 \and Yeong-Nan Yeh$^{c,}$\thanks{Partially supported by NSC 96-2115-M-006-012}}
\maketitle \vspace*{-1.2cm}\begin{center} \footnotesize $^{a}$
D$\acute{e}$pt. de math$\acute{e}$matiques, Universit$\acute{e}$ du
Qu$\acute{e}$bec $\grave{a}$ Montr$\acute{e}$al, C.P. 8888, Succ.
Centre-Ville Montr$\acute{e}$al, Canada H3C 3P8\\
$^{b,c}$Institute of Mathematics, Academia Sinica, Taipei, Taiwan
\end{center}
\date{}
 \vspace*{-0.3cm}
\thispagestyle{empty}
%--------------------------------------------------
%Abstract
%--------------------------------------------------
\begin{abstract}
The circular peak set of a permutation $\sigma$ is the set
$\{\sigma(i)\mid \sigma(i-1)<\sigma(i)>\sigma(i+1)\}$. Let
$\mathcal{P}_n$ be the set of all the subset $S\subseteq [n]$ such
that there exists a permutation $\sigma$ which has the circular set
$S$. We can make the set $\mathcal{P}_n$ into a poset
$\mathscr{P}_n$ by defining $S\preceq T$ if $S\subseteq T$ as sets.
In this paper, we prove that the poset $\mathscr{P}_n$ is a
simplicial complex on the vertex set $[3,n]$. We study the
$f$-vector, the $f$-polynomial, the reduced Euler characteristic,
the M$\ddot{o}$bius function, the $h$-vector and the $h$-polynomial
of $\mathscr{P}_n$. We also derive the zeta polynomial of
$\mathscr{P}_n$ and give the formula for the number of the chains in
$\mathscr{P}_n$. By the poset $\mathscr{P}_n$, we define two
algebras $\mathcal{A}_{\mathscr{P}_n}$ and
$\mathcal{B}_{\mathscr{P}_n}$. We consider the Hilbert polynomials
and the Hilbert series of the algebra $\mathcal{A}_{\mathscr{P}_n}$
and $\mathcal{B}_{\mathscr{P}_n}$.
\end{abstract}
%Keyword
\noindent {\bf Keyword: Circular peak; Hilbert Polynomial; Hilbert
series; Permutation; Poset; Simplicial complex; }

%main text
%Section 1
\newpage
\section{Introduction}
Throughout this paper, let $[m,n]:=\{m,m+1,\cdots,n\}$, $[n]:=[1,n]$
and $\mathfrak{S}_n$ be the set of all the permutations in the set
$[n]$. We will write permutations of $\mathfrak{S}_n$ in the form
$\sigma=(\sigma(1)\sigma(2)\cdots\sigma(n))$.

We say that a permutation $\sigma$ has a { circular descent} of
value $\sigma(i)$ if $\sigma(i)>\sigma(i+1)$ for any $i\in[n-1]$.
The { circular descent set} of a permutation $\sigma$, denoted
$CDES(\sigma)$, is the set $\{\sigma(i)\mid
\sigma(i)>\sigma(i+1)\}.$  For any $S\subseteq [n]$, define a set
$CDES_n(S)$ as $CDES_n(S)=\{\sigma\in\mathfrak{S}_n\mid
CDES(\sigma)=S\}$ and use $cdes_n(S)$ to denote the number of the
permutations in the set $CDES_n(S)$, i.e., $cdes_n(S)=|CDES_n(S)|$.
In a join work \cite{chuang}, Hungyung Zhang et al. derive the
explicit formula for $cdes_n(S)$. As a application of the main
results in \cite{chuang}, they also give the enumeration of
permutation tableaux according to their shape and generalizes the
results in \cite{Dom}. Moreover, Robert J.Clarke et al.
\cite{clarke} gave the conceptions of linear peak and cyclic peak
 and studied some new Mahonian permutation
statistics. In this paper, we are interested in the circular peaks
of permutations.  A permutation $\sigma$ has a {\it circular peak}
of value $\sigma(i)$ if $\sigma(i-1)<\sigma(i)>\sigma(i+1)$ for any
$i\in[2,n-1]$. The {\it circular peak set} of a permutation
$\sigma$, denoted $CP(\sigma)$, is the set $\{\sigma(i)\mid
\sigma(i-1)<\sigma(i)>\sigma(i+1)\}.$ For example, the circular peak
set of $\sigma=(48362517)$ is $\{5,6,8\}$. Since $\sigma$ has no
circular peaks when $n\leq 2$, we may always suppose that $n\geq 3$.
For any $S\subseteq [n]$, define a set $CP_n(S)$ as
$CP_n(S)=\{\sigma\in\mathfrak{S}_n\mid CP(\sigma)=S\}.$ Obviously,
if $\{1,2\}\subseteq S$, then $CP_n(S)=\emptyset$.
\begin{exa} \begin{eqnarray*}CP_5(\{4,5\})=\{&14253,14352,24153,34152,24351,34251,&\\
&15243,15342,25143,35142,25341,35241&\}\end{eqnarray*}
\end{exa}
Suppose that $S=\{i_1,i_2,\cdots,i_k\}$, where $i_1<i_2<\cdots<i_k$.
We prove that the necessary and sufficient conditions for
$CP_n(S)\neq \emptyset$ are $i_j\geq 2j+1$ for all $j\in [k]$.

Let $\mathcal{P}_n=\{S\mid CP_n(S)\neq\emptyset\}$. we can make the
set $\mathcal{P}_n$ into a poset $\mathscr{P}_n$ by defining
$S\preceq T$ if $S\subseteq T$ as sets. We draw the Hasse diagrams
of $\mathscr{P}_3$, $\mathscr{P}_4$ and $\mathscr{P}_5$ as follows.
\begin{center}
\includegraphics[width=2in]{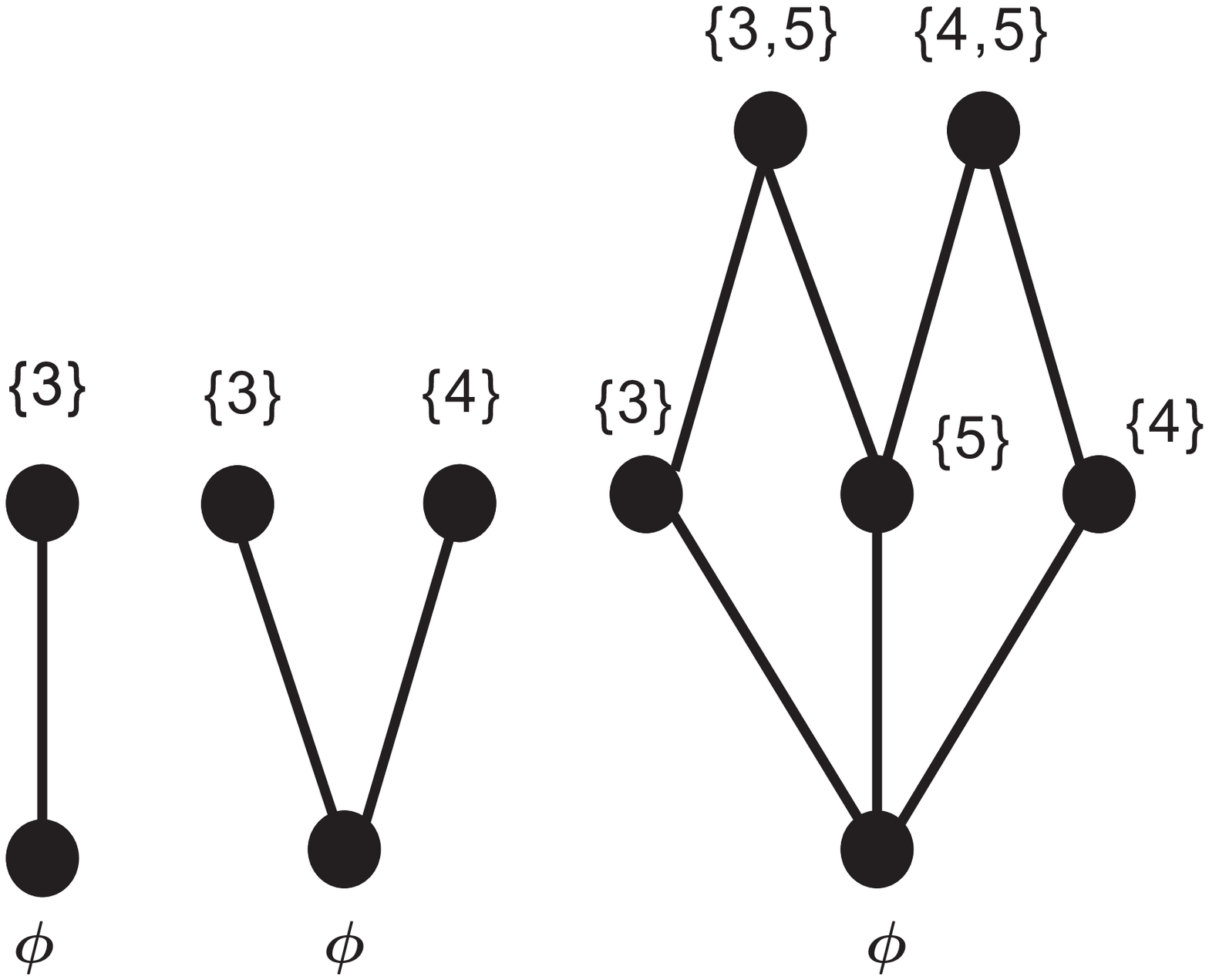}\\
 Fig.1. the Hasse
diagrams of $\mathscr{P}_3$, $\mathscr{P}_4$ and $\mathscr{P}_5$
\end{center}
There is a very interesting results: $\mathscr{P}_n$ is a simplicial
complex on the vertex set $[3,n]$. It is easy to obtain that the
dimension of the simplicial complex $\mathscr{P}_n$ is
$\lfloor\frac{1}{2}(n-1)\rfloor-1$. As we will see, the number of
the elements in $\mathscr{P}_n$ involves in the left factors of Dyck
paths of length $n-1$, counted by the ($n-1$)-th central binomial
coefficients $b_{n-1}$ (see Cori and Viennot \cite{cori}), where
$b_{n-1}={n-1\choose{\lfloor\frac{1}{2}(n-1)}\rfloor}$; the number
$p_{n,i}$ of the faces of dimension $i$ in $\mathscr{P}_n$ equals
the number of the left factors of Dyck paths from $(0,0)$ to
$(n-1.n-2i-3)$, counted by
$b_{n-1,i+1}=\frac{n-2i-2}{i+1}{n-1\choose{i}}$ \cite{cori}.

We derive the recurrence relations for the poset $\mathscr{P}_n$:
$\mathscr{P}_{n+1}\cong{\bf 2}\times\mathscr{P}_n$ if $n$ is even;
$\mathscr{P}_{n+1}\cong({\bf
2}\times\mathscr{P}_n)\setminus(\{1\}\times\mathcal{P}_{n,\lfloor\frac{1}{2}(n-1)\rfloor-1})$
if $n$ is odd, where the notation ${\bf n}$ denotes a poset formed
by the set $[n]$ with its usual order.

It is very important to obtain the $f$-vector, the $f$-polynomial
and the reduced Euler characteristic of a simplicial complex. The
integral sequence
$(p_{n,-1},p_{n,0},\cdots,p_{n,\lfloor\frac{1}{2}(n-1)\rfloor-1})$
is called the {\it f-vector} of $\mathscr{P}_n$. The {\it
$f$-polynomial} of $\mathscr{P}_n$ is defined to be the polynomial
$\mathscr{P}_n(x)=\sum\limits_{i=0}^{\lfloor\frac{1}{2}(n-1)\rfloor}p_{n,i-1}x^{\lfloor\frac{1}{2}(n-1)\rfloor-i}$.
We give the recurrence relations for the $f$-vector and  the
$f$-polynomial of $\mathscr{P}_n$. Let
$\mathscr{P}(x,y)=\sum\limits_{n\geq 3}\mathscr{P}_n(x)y^n,$ then
$\mathscr{P}(x,y)=\frac{xy^2(x+2)-xy^2C(y^2)}{x-(x+1)y^2}\frac{1+y+xy}{x+1}-y^2,$
where $C(y)=\frac{1-\sqrt{1-4y}}{2y}$. The {\it reduced Euler
characteristic} $\tilde{\chi}(\mathscr{P}_n)$ of $\mathscr{P}_n$
satisfies that $\tilde{\chi}(\mathscr{P}_n)=0$ if $n$ is odd;
$\tilde{\chi}(\mathscr{P}_n)=\frac{2(-1)^{\frac{n}{2}}}{n}{n-2\choose{\frac{1}{2}(n-2)}}$
if $n$ is even. As a corollary, the zeta polynomial of
$\mathscr{P}_n$ is
$(m-1)^{\lfloor\frac{n-1}{2}\rfloor}\mathscr{P}_n(\frac{1}{m-1})$,
which is  the number of multichains $S_{n,1}\preceq S_{n,2}\preceq
\cdots \preceq S_{n,m-1}$ in $\mathscr{P}_n$ for any $m\geq 2$.  The
number
 of the chains $S_{n,1}\prec S_{n,2}\prec\cdots \prec
S_{n,i} $ in $\mathscr{P}_n$ equals
\begin{eqnarray*}\sum{n\choose{d_1,d_2,\cdots,d_{i+1}}}\frac{2d_{i+1}-n}{n},\end{eqnarray*}
where the sum is over all $(d_1,\cdots,d_{i+1})$ such that
$\sum\limits_{k=1}^{i+1}d_k=n$, $d_1\geq 0$, $d_k\geq 1$ for all
$2\leq k\leq i$ and $d_{i+1}\geq n-\lfloor\frac{n-1}{2}\rfloor$.

We are interested in the $h$-vector and the $h$-polynomial of the
simplicial complex $\mathscr{P}_n$. We obtain the explicit formula
for $h$-vector and give the recurrence relations for the $h$-vector
and the $h$-polynomial
 of $\mathscr{P}_n$.

We fix a field $\mathbb{K}$ and let
$m={n-1\choose{\lfloor\frac{1}{2}(n-1)\rfloor}}$ and all the
elements of $\mathscr{P}_n$ be $S_{n,1},S_{n,2},\cdots,S_{n,m}$ for
any $n\geq 3$. Let $\mathcal{I}_{\mathscr{P}_n}$ be the ideal of the
polynomial ring $\mathbb{K}[x_1,\cdots,x_m]$ which is generated by
all polynomials $x_{i}x_j$ such that $S_{n,i}$ and $S_{n,j}$ are
incomparable in the poset $\mathscr{P}_n$. Let
$\mathcal{J}_{\mathscr{P}_n}$ be the ideal of the polynomial ring
$\mathbb{K}[x_1,\cdots,x_m]$ which is generated by all polynomials
$x_{i}x_j$ such that $S_{n,i}$ and $S_{n,j}$ are incomparable in
$\mathscr{P}_n$ and the polynomials $x_i^2$ for all $i\in [m]$.
Define the two algebras $\mathcal{A}_{\mathscr{P}_n}$ and
$\mathcal{B}_{\mathscr{P}_n}$ as the quotients
$\mathcal{A}_{\mathscr{P}_n}=\mathbb{K}[x_1,\cdots,x_m]/
\mathcal{I}_{\mathscr{P}_n}$ and
$\mathcal{B}_{\mathscr{P}_n}=\mathbb{K}[x_1,\cdots,x_m]/
\mathcal{J}_{\mathscr{P}_n}$, respectively. We study the Hilbert
polynomials and the Hilbert series of the algebras
$\mathcal{A}_{\mathscr{P}_n}$ and $\mathcal{B}_{\mathscr{P}_n}$.

This paper is organized as follows. In Section $2$, we give the
necessary and sufficient conditions for $CP_{n}(S)\neq \emptyset$.
In Section $3$, we prove that the poset $\mathscr{P}_n$ is a
simplicial complex and study its properties. In Section $4$, we
define the two algebras $\mathcal{A}_{\mathscr{P}_n}$ and
$\mathcal{B}_{\mathscr{P}_n}$ and consider their Hilbert polynomial
and Hilbert series.

\section{The Necessary and Sufficient Conditions for $CP_n(S)\neq \emptyset$ }
In this section, we will give the necessary and sufficient
conditions for $CP_n(S)\neq \emptyset$ for any $n\geq 3$ and
$S\subseteq [n]$.

\begin{thm}\label{theoremnecessarysufficientconditions} Suppose that $n$ is an integer with $n\geq 3$. Let
$S=\{i_1,i_2,\cdots,i_k\}$ be a subset of $[n]$, where
$i_1<i_2<\cdots<i_k$. Then the necessary and sufficient conditions
for $CP_n(S)\neq \emptyset$  are $i_j\geq 2j+1$ for all $j\in [k]$.
\end{thm}
\begin{proof} First, we suppose that $CP_n(S)\neq \emptyset$ and
let $\sigma\in CP_n(S)$. For any $j\in[k]$, since all the numbers
$i_1,i_2,\cdots,i_j$ are peaks of $\sigma$, we have $i_j-j\geq j+1$,
hence, $i_j\geq 2j+1$.

Conversely, when $i_j\geq 2j+1$ for all $j\in [k]$, let
$T=[i_k]\setminus S$, then the set $T$ has at least $k+1$ elements.
So, suppose that $T=\{a_1,a_2,\cdots,a_m\}$ with $a_1<a_2<\cdots
<a_m$. We consider $\sigma=a_1i_1a_2i_2\cdots a_ki_ka_{k+1}\cdots
a_m(i_k+1)\cdots n$. Obviously, $CP(\sigma)\subseteq S$. If
$CP(\sigma)\neq S$, then there exists a minimal $j\in[k]$ such that
$i_j$ is not a circular peak of $\sigma$. So, $a_{j+1}>i_j$. This
implies that $i_j=2j$, a contradiction. Hence, $CP(\sigma)=S$ and
$CP_n(S)\neq\emptyset$. \end{proof}

\begin{cor}\label{corollaryupperboundofs}Suppose that $n$ is an integer with $n\geq 3$ and
$S\subseteq [n]$. If $CP_n(S)\neq\emptyset$, then $|S|\leq
\lfloor\frac{n-1}{2}\rfloor$.
\end{cor}
\begin{proof} Suppose that $S=\{i_1,i_2,\cdots,i_k\}$ with
$i_1<i_2<\cdots<i_k$. Since $CP_n(S)\neq\emptyset$, Theorem
\ref{theoremnecessarysufficientconditions} tells us that $i_k\geq
2k+1$. $i_k\leq n$ implies that $k\leq
\lfloor\frac{n-1}{2}\rfloor$.\end{proof}

\begin{cor}\label{corScupn+1} Let $n\geq 3$ and $S\subseteq [n]$. Suppose that
$CP_n(S)\neq\emptyset$. Then when $|S|<\lfloor\frac{n-1}{2}\rfloor$,
we have $CP_{n+1}(S\cup\{n+1\})\neq\emptyset$; when
$|S|=\lfloor\frac{n-1}{2}\rfloor$, we have
$CP_{n+1}(S\cup\{n+1\})\neq\emptyset$ if $n$ is even; otherwise,
$CP_{n+1}(S\cup\{n+1\})=\emptyset$.
\end{cor}
\begin{proof} Suppose that $k=|S|$. $k<\lfloor\frac{n-1}{2}\rfloor$ implies
that $2(k+1)+1\leq 2\lfloor\frac{n-1}{2}\rfloor+1<n+1$. So,
$CP_{n+1}(S\cup\{n+1\})\neq\emptyset$ when
$|S|<\lfloor\frac{n-1}{2}\rfloor$. For the case
$k=\lfloor\frac{n-1}{2}\rfloor$, we have
$$2(k+1)+1=\left\{\begin{array}{lll}n+1&\text{if}&n \text{ is even}\\
n+2&\text{if}&n\text{ is odd}\end{array}\right..$$ By Theorem
\ref{theoremnecessarysufficientconditions},
$CP_{n+1}(S\cup\{n+1\})\neq\emptyset$ if $n$ is even; otherwise,
$CP_{n+1}(S\cup\{n+1\})=\emptyset$.\end{proof}

\section{The Simplicial Complex $\mathscr{P}_n$}
In this section, we will prove that the poset $\mathscr{P}_n$ is a
simplicial complex on the vertex set $[3,n]$ for any $n\geq 3$, and
then study the properties of the Simplicial Complex $\mathscr{P}_n$.
Following \cite{Stanley}, we define a {\it simplicial complex}
$\Delta$ on a vertex set
$V$ as a collection of subsets of $V$ satisfying:\\
(1) If $x\in V$, then $\{x\}\in \Delta$, and\\
(2) if $S\in \Delta$ and $T\subseteq S$, then $T\in\Delta$.

\begin{thm}Let $n\geq 3$. Then $\mathscr{P}_n$ is a simplicial complex of the set
$[3,n]$ and has the dimension $\lfloor\frac{n-1}{2}\rfloor-1$.
\end{thm}
\begin{proof} Obviously, $\emptyset\in\mathscr{P}_n$. For any $3\leq x\leq
n$, Theorem \ref{theoremnecessarysufficientconditions} implies that
$\{x\}\in\mathscr{P}_n$.   Noting that if $CP_n(T)=\emptyset$ then
$CP_n(S)=\emptyset$ for any $S\supseteq T$, we have if $S\in
\mathcal{P}_n$ and $T\subseteq S$ then $T\in\mathscr{P}_n$. Hence,
$\mathscr{P}_n$ is a simplicial complex of the set
$[3,n]$.\end{proof}

An element $S\in\mathscr{P}_n$ is called a {\it face} of
$\mathscr{P}_n$, and the {\it dimension} of $S$ is defined to be
$|S|-1$. In particular, the void set $\emptyset$ is always a face of
$\mathscr{P}_n$ of dimension $-1$. Also define the dimension of
$\mathscr{P}_n$ by
$dim\mathscr{P}_n=\max\limits_{S\in\mathscr{P}_n}(dim S).$

\begin{thm}The simplicial complex $\mathscr{P}_n$ has the dimension $\lfloor\frac{n-1}{2}\rfloor-1$.
\end{thm}
\begin{proof} Taking $S=\{3,5,\cdots,2\lfloor\frac{n-1}{2}\rfloor+1\}$, by
Theorem \ref{theoremnecessarysufficientconditions}, we have
$S\in\mathscr{P}_n$. From Corollary \ref{corollaryupperboundofs} it
follows  that the dimension of $\mathscr{P}_n$ is
$\lfloor\frac{n-1}{2}\rfloor-1$.\end{proof}

There are very close relations between the number of the elements of
$\mathscr{P}_n$ and the number of left factor of Dyck path of length
$n$. An $n$-\emph{Dyck path} is a lattice path in the first quadrant
starting at $(0,0)$ and ending at $(2n,0)$ with only two kinds of
steps---\emph{rise step: $U=(1,1)$ }and \emph{fall step:
$D=(1,-1)$}. We can also consider an $n$-Dyck path $P$ as a word of
$2n$ letters using only $U$ and $D$. Let $L=w_1w_2\cdots w_n$ be a
word, where $w_j\in\{U,D\}$ and $n\geq 0$. If there is another word
$R$ which consists of $U$ and $D$ such that $LR$ forms a Dyck path,
then $L$ is called an $n$-left factor of Dyck path. Let
$\mathcal{L}_n$ denote the set of all $n$-left factor of Dyck paths.
It is well known that $|\mathcal{L}_n|$, the cardinality of
$\mathcal{L}_n$, equals the $n$th Central binomial number
$b_n={n\choose \lfloor \frac{n}{2}\rfloor}$. In the following lemma,
we give a bijection $\phi$ from the sets $\mathcal{P}_n$ to
$\mathcal{L}_{n-1}$.

\begin{lem}\label{theorembijectionleftor}There is a bijection $\phi$ between the sets $\mathcal{P}_n$ and
$\mathcal{L}_{n-1}$ for any $n\geq 3$. Furthermore, the number of
the elements in $\mathscr{P}_n$ is ${n-1\choose \lfloor
\frac{n-1}{2}\rfloor}$.
\end{lem}
\begin{proof} For any $S\in\mathcal{P}_n$, we construct a word
$\phi(S)=w_1w_2\cdots w_{n-1}$ as follows:
$$w_i=\left\{\begin{array}{lll}
D&\text{if}&i+1\in S\\
U&\text{if}&i+1\notin S
\end{array}\right.$$ for any $i\in[n-1]$. Theorem \ref{theoremnecessarysufficientconditions}
implies that $\phi(S)$ is a $(n-1)$-left factor of Dyck path.
Conversely, for any a $n$-left factor of Dyck path $w_1w_2\cdots
w_{n-1}$, let $S=\{i+1\mid w_i=D\}$, then $CP_n(S)\neq\emptyset$.
Hence, the mapping $\phi$ is a bijection. Note that the number of
$(n-1)$-left factor of Dyck path is ${n-1\choose \lfloor
\frac{n-1}{2}\rfloor}$. Hence, $|\mathcal{P}_n|={n-1\choose \lfloor
\frac{n-1}{2}\rfloor}$. \end{proof}

Now, we are in a position to obtain the number $p_{n,i}$ of the
faces of dimension $i$ in $\mathscr{P}_n$. For any $i\geq 0$, let
$\mathcal{L}_{n,i}$ denote the set of all $n$-left factor of Dyck
paths from $(0,0)$ to $(n,n-2i)$. Define a set $\mathcal{P}_{n,i}$
as the set of all the faces of dimension $i$ in $\mathscr{P}_n$,
i.e., $\mathcal{P}_{n,i}=\{S\in\mathcal{P}_n\mid |S|=i+1\}$ for any
$-1\leq i\leq \lfloor\frac{n-1}{2}\rfloor-1$. Clearly,
$p_{n,i}=|\mathcal{P}_{n,i}|$.

\begin{cor} Let $n\geq 3$. There is a bijection between the sets $\mathcal{P}_{n,i}$ and
$\mathcal{L}_{n-1,i+1}$ for any $-1\leq i\leq \lfloor
\frac{n-1}{2}\rfloor-1$. Furthermore, we have
$$p_{n,i}=\left\{\begin{array}{lll}1&if&i=-1\\
\frac{n-2i-2}{i+1}{n-1\choose{i}}&if&0\leq i\leq \lfloor
\frac{n-1}{2}\rfloor-1 \end{array}\right.$$
\end{cor}
\begin{proof} We just consider the case with $i\geq 0$. For any
$S\in\mathcal{P}_{n,i}$, $|S|=i+1$ implies that the number of letter
$D$ in the word $\phi(S)$ is  $i+1$. Hence, $\phi(S)$ is a left
factor of Dyck path from $(0,0)$ to $(n-1,n-2i-3)$. So,
$\phi(S)\in\mathcal{L}_{n-1,i+1}$. \end{proof}

If $P$ and $Q$ are posets, then the {\it direct product} of $P$ and
$Q$ is the poset $P\times Q$ on the set $\{(x,y)\mid x\in P\text{
and }y\in Q\}$ such that $(x,y)\leq (x',y')$ in $P\times Q$ if
$x\leq x'$ in $P$ and $y\leq y'$ in $Q$. Recall that the poset $\bf
n$ is formed by the set $[n]$ with its usual order. By Corollary
\ref{corScupn+1}, we obtain a method for constructing the poset
$\mathscr{P}_{n+1}$ from $\mathscr{P}_{n}$.

\begin{thm}\label{lemmaposetdiagram}$\mathscr{P}_{n+1}\cong{\bf
2}\times\mathscr{P}_n$ if $n$ is even; $\mathscr{P}_{n+1}\cong({\bf
2}\times\mathscr{P}_n)\setminus(\{1\}\times\mathcal{P}_{n,\lfloor\frac{n-1}{2}\rfloor-1})$
if $n$ is odd.
\end{thm}

By Theorem \ref{lemmaposetdiagram}, it is easy for us to obtain the
M$\ddot{o}$bius function, the recurrence relations for the
$f$-vector and  the $f$-polynomial of the poset $\mathscr{P}_n$.

\begin{cor} Let $\mu_{\mathscr{P}_n}$ be the M$\ddot{o}$bius
function of the poset $\mathscr{P}_n$. Then
\begin{eqnarray*}\mu_{\mathscr{P}_n}(S,T)=(-1)^{|T|-|S|}\end{eqnarray*}for
any $S\preceq T$ in $\mathscr{P}_n$.
\end{cor}
\begin{proof} Obviously, $\mu_{\mathscr{P}_3}(\emptyset,\{3\})=-1$. By induction for $n$, we assume that
$\mu_{\mathscr{P}_n}(S,T)=(-1)^{|T|-|S|}$ for any $S\preceq T$ in
$\mathscr{P}_n$. Lemma \ref{lemmaposetdiagram} tells us that
$\mathscr{P}_{n+1}\cong{\bf 2}\times\mathscr{P}_n$ if $n$ is even;
$\mathscr{P}_{n+1}\cong({\bf
2}\times\mathscr{P}_n)\setminus(\{1\}\times\mathcal{P}_{n,\lfloor\frac{n-1}{2}\rfloor})$
if $n$ is odd. We conclude from the product theorem that
$$\mu_{\mathscr{P}_{n+1}}(S,T)=\left\{\begin{array}{lll}
\mu_{\mathscr{P}_n}(S\setminus
\{n+1\},T\setminus\{n+1\})&\text{if}&n+1\in S\cap
T\\
\mu_{\mathscr{P}_n}(S,T)&\text{if}&n+1\notin S\cup
T\\
-\mu_{\mathscr{P}_n}(S,T\setminus\{n+1\})&\text{if}&n+1\notin S
\text{ and } n+1\in T
\end{array}\right.$$ for any $S\prec T$. Note that $\mu_{\mathscr{P}_n}(S\setminus
\{n+1\},T\setminus\{n+1\})=(-1)^{|T|-1-(|S|-1)}=(-1)^{|T|-|S|}$ and
$-\mu_{\mathscr{P}_n}(S,T\setminus\{n+1\})=-(-1)^{|T|-1-|S|}=(-1)^{|T|-|S|}$.
Hence, $\mu_{\mathscr{P}_{n+1}}(S,T)=(-1)^{|T|-|S|}$.\end{proof}

\begin{cor}\label{corfrec} Let $n\geq 3$. The sequence $p_{n,i}$ satisfies the
following recurrence relation: when $n$ is even
,$$p_{n+1,i}=\left\{\begin{array}{lll}
p_{n,i}&\text{if}&i=-1\\
p_{n,i-1}+p_{n,i}&\text{if}&i=0,1,\cdots,\frac{n}{2}-2\\
p_{n,i-1}&\text{if}&i=\frac{n}{2}-1\\
\end{array}\right.;$$
when $n$ is odd ,$$p_{n+1,i}=\left\{\begin{array}{lll}
p_{n,i}&\text{if}&i=-1\\
p_{n,i-1}+p_{n,i}&\text{if}&i=0,1,\cdots,\frac{n-3}{2}
\end{array}\right.$$ with initial
conditions $(p_{3,-1},p_{3,0})=(1,1)$.
\end{cor}
\begin{proof} First, we consider the case with $n$ even. It is easy
to check that $p_{n+1,-1}=p_{n,-1}=1$. For any
$S\in\mathcal{P}_{n+1,\frac{1}{2}n-1}$,  Corollary \ref{corScupn+1}
tells us that $n+1\in S$. Note that
$S\in\mathcal{P}_{n+1,\frac{1}{2}n-1}$ if and only if
$S\setminus\{n+1\}\in\mathcal{P}_{n,\frac{1}{2}n-2}$. Hence,
$p_{n+1,\frac{1}{2}n-1}=p_{n,\frac{1}{2}n-2}$. Let $0\leq i\leq
\frac{1}{2}n-2$,  obviously, $\mathcal{P}_{n,i}\subseteq
\mathcal{P}_{n+1,i}$. For any $S\in\mathcal{P}_{n+1,i}$ with $n+1\in
S$, $S\setminus \{n+1\}$ can be viewed as a element of
$\mathcal{P}_{n,i-1}$. Conversely, for any
$S\in\mathcal{P}_{n,i-1}$, Corollary \ref{corScupn+1} implies that
$S\cup\{n+1\}\in\mathcal{P}_{n+1,i}$. Hence,
$p_{n+1,i}=p_{n,i-1}+p_{n,i}$. Similarly, we can consider the case
with $n$ odd. \end{proof}

\begin{thm}\label{theoremfpolynomial} The $f$-polynomial $\mathscr{P}_n(x)$ of the simplicial complex
$\mathscr{P}_n$ satisfy the following recurrence relation: when $n$
is even, then
\begin{eqnarray*}\mathscr{P}_{n+1}(x)
=(1+x)\mathscr{P}_{n}(x),\end{eqnarray*} and when $n$ is odd,  then
\begin{eqnarray*}x\mathscr{P}_{n+1}(x)
=(1+x)\mathscr{P}_{n}(x)-\frac{2}{n+1}{n-1\choose{\frac{n-1}{2}}}.\end{eqnarray*}
for any $n\geq 3$, with initial condition $\mathscr{P}_3(x)=x+1$.

Let $\mathscr{P}(x,y)=\sum\limits_{n\geq 3}\mathscr{P}_n(x)y^n.$
Then
$$\mathscr{P}(x,y)=\frac{xy^2(x+2)-xy^2C(y^2)}{x-(x+1)y^2}\frac{1+y+xy}{x+1}-y^2,$$
where $C(y)=\frac{1-\sqrt{1-4y}}{2y}$.
\end{thm}
\begin{proof} Obviously, $\mathscr{P}_3(x)=x+1$. When $n$ is odd, we
suppose that $n=2i+1$ with $i\geq 1$. Corollary \ref{corfrec}
implies that \begin{eqnarray*}x\mathscr{P}_{2i+2}(x)
=(1+x)\mathscr{P}_{2i+1}(x)-\frac{1}{(i+1)}{2i\choose{i}}.\end{eqnarray*}
Similarly, when $n$ is even, we suppose that $n=2i$ with $i\geq 2$.
By corollary \ref{corfrec}, we have
\begin{eqnarray*}\mathscr{P}_{2i+1}(x)
=(1+x)\mathscr{P}_{2i}(x).\end{eqnarray*}

Let $\mathscr{P}_{odd}(x,y)=\sum\limits_{i\geq
1}\mathscr{P}_{2i+1}(x)y^{2i+1}$ and
$\mathscr{P}_{even}(x,y)=\sum\limits_{m\geq
2}\mathscr{P}_{2i}(x)y^{2i}$, then
$\mathscr{P}_{odd}(x,y)=(x+1)y^3+(x+1)y\mathscr{P}_{even}(x,y)$ and
$\mathscr{P}(x,y)=\mathscr{P}_{odd}(x,y)+\mathscr{P}_{even}(x,y)$.
It is easy to check that
$x\mathscr{P}_{2i+3}(x)=(1+x)^2\mathscr{P}_{2i+1}(x)-\frac{1}{i+1}{2i\choose{i}}x$.
So, $\mathscr{P}_{odd}$ satisfies the following equation
\begin{eqnarray*}x\mathscr{P}_{odd}(x,y)=(x+1)^2y^2\mathscr{P}_{odd}(x,y)+x(x+2)y^3-xy^3C(y^2),\end{eqnarray*}where
$C(y)=\frac{1-\sqrt{1-4y}}{2y}$. Equivalently,
$$\mathscr{P}_{odd}(x,y)=\frac{xy^3(x+2)-xy^3C(y^2)}{x-(x+1)y^2}.$$
Hence,
$$\mathscr{P}(x,y)=\frac{xy^2(x+2)-xy^2C(y^2)}{x-(x+1)y^2}\frac{1+y+xy}{x+1}-y^2.$$
\end{proof}

Define the {\it reduced Euler characteristic} of $\mathscr{P}_n$ by
$\tilde{\chi}(\mathscr{P}_n)=\sum\limits_{i=0}^{\lfloor\frac{1}{2}(n-1)\rfloor}(-1)^{i-1}p_{n,i-1}$.

\begin{cor} For any $n\geq 3$, $\tilde{\chi}(\mathscr{P}_n)=\left\{\begin{array}{lll}0&\text{\it if}&n~\text{\it is odd}\\\frac{2(-1)^{\frac{n}{2}}}{n}{n-2\choose{\frac{1}{2}(n-2)}}&\text{\it if}&n~\text{\it is even}\end{array}\right.$
\end{cor}
\begin{proof} Clearly, $\mathscr{P}_{3}(-1)=0$. Theorem
\ref{theoremfpolynomial} tells us that
$$\mathscr{P}_{n+1}(-1)=\left\{\begin{array}{lll}0&\text{if}&n\text{ is even}\\\frac{2}{n+1}{n-1\choose{\frac{1}{2}(n-1)}}&\text{if}&n\text{ is odd}\end{array}\right.$$ for any $n\geq 4$. Since
$\tilde{\chi}(\mathscr{P}_n)=(-1)^{\lfloor\frac{n-1}{2}\rfloor-1}\mathscr{P}_n(-1)$,
we have
$$\tilde{\chi}(\mathscr{P}_n)=\left\{\begin{array}{lll}0&\text{if}&n~\text{\it is odd}\\\frac{2(-1)^{\frac{n}{2}-2}}{n}{n-2\choose{\frac{1}{2}(n-2)}}&\text{if}&n~\text{\it is even}\end{array}\right.$$  \end{proof}

Let $P$ be a finite post. Define $Z(P,i)$ to be the number of
multichains $x_1\leq x_2\leq \cdots \leq x_{i-1}$ in $P$ for any
$i\geq 2$. $Z(P,i)$ is called the {\it zeta polynomial} of $P$. We
state Proposition $3.11.1a$ and Proposition  $3.14.2$ in
\cite{Stanley} as the following lemma.
\begin{lem}\label{rankfunction} Suppose that $P$ is a poset.
(1) Let $d_i$ be the number of chains $x_1<x_2<\cdots <x_{i-1}$ in
$P$. Then $Z(P,i)=\sum\limits_{j\geq 2}d_j{i-2\choose{j-2}}$.

(2) If $P$ is simiplicial and graded, then $Z(P,x+1)$ is the
rank-generating function of $P$.
\end{lem}

\begin{cor}\label{zetefunction} Let $n\geq 3$ and $m\geq 2$. Then
$Z(\mathscr{P}_n,i)=(i-1)^{\lfloor\frac{n-1}{2}\rfloor}\mathscr{P}_n(\frac{1}{i-1})$
for any $i\geq 2$. Furthermore, $Z(\mathscr{P}_n,i)$ satisfies the
recurrence relations:
$Z(\mathscr{P}_{n+1},i)=iZ(\mathscr{P}_n,i)-\varepsilon(n)\frac{2(i-1)^{\frac{1}{2}(n+1)}}{n+1}{n-1\choose{\frac{1}{2}(n-1)}}$,
where $\varepsilon(n)=0$ if $n$ is even; $\varepsilon(n)=1$
otherwise, with initial condition $Z(\mathscr{P}_3,i)=i$.\end{cor}
\begin{proof} Let $\mathscr{P}_n(x)$ be the $f$-polynomial of
$\mathscr{P}_n$, then the rank-generating function of
$\mathscr{P}_n$ is
$x^{\lfloor\frac{1}{2}(n-1)\rfloor}\mathscr{P}_n(\frac{1}{x})$.
Lemma \ref{rankfunction}(2) implies that
$Z(\mathscr{P}_n,i)=(i-1)^{\lfloor\frac{n-1}{2}\rfloor}\mathscr{P}_n(\frac{1}{i-1})$.
The recurrence relations for $Z(\mathscr{P}_n,i)$ follows Theorem
\ref{theoremfpolynomial}. \end{proof}

Let $d_{\mathscr{P}_{n},i}$ be the number of the chains
$S_{n,1}\prec S_{n,2}\prec\cdots \prec S_{n,i} $ of $\mathscr{P}_n$.

\begin{thm}\label{thmchain} For any  $i\geq 1$,
\begin{eqnarray*}d_{\mathscr{P}_{n},i}=\sum{n\choose{d_1,d_2,\cdots,d_{i+1}}}\frac{2d_{i+1}-n}{n},\end{eqnarray*}
where the sum is over all $(d_1,\cdots,d_{i+1})$ such that
$\sum\limits_{k=1}^{i+1}d_k=n$, $d_1\geq 0$, $d_k\geq 1$ for all
$2\leq k\leq i$ and $d_{i+1}\geq n-\lfloor\frac{n-1}{2}\rfloor$.
\end{thm}
\begin{proof} Let  $i\geq 1$ and $S_{n,1}\prec S_{n,2}\prec\cdots
\prec S_{n,i} $ be a chain of  $\mathscr{P}_n$. Suppose that
$|S_{n,k}|=j_k$ for any $k\in[i]$, then $0\leq j_1<j_2<\cdots
<j_{i}\leq \lfloor\frac{n-1}{2}\rfloor$. There are $p_{n,j_{i}-1}$
ways to obtain the set $S_{n,i}$. Given $S_{n,k}$ with  $k\geq 2$,
there are ${j_k\choose{j_{k-1}}}$ ways to form the subset
$S_{n,k-1}\subseteq S_{n,k}$. Hence,
\begin{eqnarray*}d_{\mathscr{P}_{n},i}&=&\sum\limits_{0=j_0\leq j_1<j_2<\cdots < j_{i}\leq \lfloor\frac{n-1}{2}\rfloor}
\prod\limits_{k=0}^{i-1}{j_{k+1}\choose{j_k}}p_{n,{j_i-1}}\\
&=&\sum{n\choose{d_1,d_2,\cdots,d_{i+1}}}\frac{2d_{i+1}-n}{n},\end{eqnarray*}
where the sum is over all $(d_1,\cdots,d_{i+1})$ such that
$\sum\limits_{k=1}^{i+1}d_k=n$, $d_1\geq 0$, $d_k\geq 1$ for all
$2\leq k\leq i$ and $d_{i+1}\geq n-\lfloor\frac{n-1}{2}\rfloor$.
\end{proof}

\begin{cor}For any $n\geq 3$,
$$\mathscr{P}_n(x)=\sum\limits_{i=2}^{\lfloor\frac{n-1}{2}\rfloor+2}\frac{x^{\lfloor\frac{n-1}{2}\rfloor+2-i}}{(i-2)!}\prod\limits_{j=1}^{i-2}(1-jx)\sum{n\choose{d_1,d_2,\cdots,d_{i}}}\frac{2d_{i}-n}{n}$$where the second sum is over all $(d_1,\cdots,d_{i})$ such that
$\sum\limits_{k=1}^{i}d_k=n$, $d_1\geq 0$, $d_k\geq 1$ for all
$2\leq k\leq i-1$ and $d_{i}\geq n-\lfloor\frac{n-1}{2}\rfloor$.
\end{cor}
\begin{proof} Lemma \ref{rankfunction}(1) implies that
$Z(\mathscr{P}_n,i)=\sum\limits_{j=2}^{\lfloor\frac{n-1}{2}\rfloor+2}
d_{\mathscr{P}_n,j-1}{i-2\choose{j-2}}$. By Corollary
\ref{zetefunction}, we have
$$\mathscr{P}_n\left(\frac{1}{i-1}\right)=\left(\frac{1}{i-1}\right)^{\lfloor\frac{n-1}{2}\rfloor}\sum\limits_{j=2}^{\lfloor\frac{n-1}{2}\rfloor+2}
d_{\mathscr{P}_n,j-1}{i-2\choose{j-2}}$$ for any $i\geq 2$. Note
that
$$x^{\lfloor\frac{n-1}{2}\rfloor}\sum\limits_{j=2}^{\lfloor\frac{n-1}{2}\rfloor+2}d_{\mathscr{P}_n,j-1}{\frac{1}{x}-1\choose{j-2}}=\sum\limits_{j=2}^{\lfloor\frac{n-1}{2}\rfloor+2}\frac{x^{\lfloor\frac{n-1}{2}\rfloor+2-j}}{(j-2)!}\prod\limits_{k=1}^{j-2}(1-kx)d_{\mathscr{P}_n,j-1}$$
is a polynomial. Hence,
$\mathscr{P}_n(x)=\sum\limits_{j=2}^{\lfloor\frac{n-1}{2}\rfloor+2}\frac{x^{\lfloor\frac{n-1}{2}\rfloor+2-j}}{(j-2)!}\prod\limits_{k=1}^{j-2}(1-kx)d_{\mathscr{P}_n,j-1}$.
\end{proof}

Let
$\mathscr{H}_n(x)=\mathscr{P}_n(x-1)=\sum\limits_{i=0}^{\lfloor\frac{1}{2}(n-1)\rfloor}h_{n,i}x^{\lfloor\frac{1}{2}(n-1)\rfloor-i}$,
then $\mathscr{H}_n(x)$ is called {\it $h$-polynomial} of
$\mathscr{P}_n$ and the sequence
$(h_{n,0},h_{n,1},\cdots,h_{n,\lfloor\frac{1}{2}(n-1)\rfloor})$ {\it
$h$-vector} of $\mathscr{P}_n$.

\begin{cor}\label{corhpolynomial}  The $h$-polynomial $\mathscr{H}_n(x)$ of the simplicial complex
$\mathscr{P}_n$  satisfies the  recurrence relation: when $n$ is
even,
\begin{eqnarray*}\mathscr{H}_{n+1}(x)
=x\mathscr{H}_{n}(x),\end{eqnarray*} and when $n$ is odd,
\begin{eqnarray*}(x-1)\mathscr{H}_{n+1}(x)
=x\mathscr{H}_{n}(x)-\frac{2}{(n+1)}{n-1\choose{\frac{n-1}{2}}},\end{eqnarray*}
for any $n\geq 3$, with initial condition $\mathscr{H}_3(x)=x$.

Let $\mathscr{H}(x,y)=\sum\limits_{n\geq 3}\mathscr{H}_n(x)y^n$,
then
\begin{eqnarray*}\mathscr{H}(x,y)=\frac{\left[(x^2-1)y^2-(x-1)y^2C(y^2)\right](1+xy)}{x(x-1-xy^2)}-y^2.\end{eqnarray*}

Furthermore, let
$(h_{n,0},h_{n,1},\cdots,h_{n,\lfloor\frac{n-1}{2}\rfloor})$ be the
$h$-vector
 of
$\mathscr{P}_n$, then
\begin{eqnarray*}h_{n,i}=\frac{\lfloor\frac{n}{2}\rfloor-i}{\lfloor\frac{n}{2}\rfloor+i}{\lfloor\frac{n}{2}\rfloor+i\choose{\lfloor\frac{n}{2}\rfloor}}.\end{eqnarray*}
\end{cor}
\begin{proof} Since $\mathscr{H}_n(x)=\mathscr{P}_n(x-1)$, by Theorem
\ref{theoremfpolynomial}, we easily obtain that if $n$ is even, then
\begin{eqnarray*}\mathscr{H}_{n+1}(x)
=x\mathscr{H}_{n}(x),\end{eqnarray*} and if $n$ is odd,  then
\begin{eqnarray*}(x-1)\mathscr{H}_{n+1}(x)
=x\mathscr{H}_{n}(x)-\frac{2}{n+1}{n-1\choose{\frac{n-1}{2}}},\end{eqnarray*}
for any $n\geq 3$, with initial condition $\mathscr{H}_3(x)=x$.

Since $\mathscr{H}(x,y)=\mathcal{P}(x-1,y)$, we have
\begin{eqnarray*}\mathscr{H}(x,y)=\frac{\left[(x^2-1)y^2-(x-1)y^2C(y^2)\right](1+xy)}{x(x-1-xy^2)}-y^2.\end{eqnarray*}
\end{proof}

\begin{cor} Let the sequence
$(h_{n,0},h_{n,1},\cdots,h_{n,\lfloor\frac{1}{2}(n-1)\rfloor})$ be
$h$-vector of $\mathscr{P}_n$. Then the sequence $h_{n,i}$ satisfies
the following recurrence relation:
$$\left\{\begin{array}{lllll}
h_{n+1,0}&=&h_{n,0}\\
h_{n+1,i}&=&h_{n,i}+\varepsilon(n)h_{n+1,i-1}&\text{if}&1\leq i\leq \lfloor\frac{n}{2}\rfloor-1\\
h_{n+1,\lfloor\frac{n}{2}\rfloor}&=&\varepsilon(n)c_{\lfloor\frac{n}{2}\rfloor}\end{array}\right.$$
where $c_m=\frac{1}{m+1}{2m\choose{m}}$ and $\varepsilon(n)=0$ if
$n$ is even; otherwise, $\varepsilon(n)=1$, with initial conditions
$(h_{3,0},h_{3,1})=(1,0)$. Equivalently,
\begin{eqnarray*}h_{n,i}=\frac{\lfloor\frac{n}{2}\rfloor-i}{\lfloor\frac{n}{2}\rfloor+i}{\lfloor\frac{n}{2}\rfloor+i\choose{\lfloor\frac{n}{2}\rfloor}}.\end{eqnarray*}
\end{cor}
\begin{proof} By comparing with the coefficient in Corollary
\ref{corhpolynomial}, we can obtain the desired recurrence
relations. Consider
$t_{n,i}=\frac{\lfloor\frac{n}{2}\rfloor-i}{\lfloor\frac{n}{2}\rfloor+i}{\lfloor\frac{n}{2}\rfloor+i\choose{\lfloor\frac{n}{2}\rfloor}}$.
Note  that $t_{n,i}$ satisfies the above recurrence relations as
well. Hence,
\begin{eqnarray*}h_{n,i}=t_{n,i}=\frac{\lfloor\frac{n}{2}\rfloor-i}{\lfloor\frac{n}{2}\rfloor+i}{\lfloor\frac{n}{2}\rfloor+i\choose{\lfloor\frac{n}{2}\rfloor}}.\end{eqnarray*}
\end{proof}
\begin{rem} Let $n\geq 3$. The number of left factor of Dyck path from $(0,0)$ to
$(\lfloor\frac{n}{2}\rfloor+i-1,\lfloor\frac{n}{2}\rfloor-i-1)$
equals
$\frac{\lfloor\frac{n}{2}\rfloor-i}{\lfloor\frac{n}{2}\rfloor+i}{\lfloor\frac{n}{2}\rfloor+i\choose{\lfloor\frac{n}{2}\rfloor}}$
for any $0\leq i\leq \lfloor\frac{n-1}{2}\rfloor$.\end{rem}

\section{The Algebras $\mathcal{A}_{\mathscr{P}_n}$ and
$\mathcal{B}_{\mathscr{P}_n}$ from the Poset $\mathscr{P}_n$} In
this section,  we will consider the properties of the algebras
$\mathcal{A}_{\mathscr{P}_n}$ and $\mathcal{B}_{\mathscr{P}_n}$.

Let $m={n-1\choose{\lfloor\frac{1}{2}(n-1)\rfloor}}$. We list all
the elements of $\mathscr{P}_n$ as $S_{n,1},S_{n,2},\cdots,S_{n,m}$.
For any a sequence $S_{n,j_1}, S_{n,j_2},\cdots, S_{n,j_i}$ in the
poset $\mathscr{P}_n$, let ${\bf r}(S_{n,j_1}, S_{n,j_2},\cdots,
S_{n,j_i})=(r_1,\cdots,r_m)$ be a vector such that $r_j=|\{k\mid
S_{n,j_k}=S_{n,j}\}|$. Furthermore, we can obtain a monomial
$$m(S_{n,j_1}, S_{n,j_2},\cdots, S_{n,j_i})=x_1^{r_1}x_2^{r_2}\cdots
x_m^{r_m}.$$  A sequence $S_{n,j_1}, S_{n,j_2},\cdots, S_{n,j_i}$ of
elements in the poset $\mathscr{P}_n$ forms a multichain if and only
if the monomial $m(S_{n,j_1}, S_{n,j_2},\cdots, S_{n,j_i})$ is
nonvanishing in the algebra $\mathcal{A}_{\mathscr{P}_n}$.

For a monomial ideal $\mathcal{I}$, the set of all monomials that do
not belong to $\mathcal{I}$ is a basis of the quotient of the
polynomial ring modulo $\mathcal{I}$, called the standard monomial
basis. Thus the monomials $m(S_{n,j_1}, S_{n,j_2},\cdots,
S_{n,j_i})$, where $(S_{n,j_1}, S_{n,j_2},\cdots, S_{n,j_i})$ ranges
over the multichains of the poset $\mathscr{P}_{n}$, form the
standard monomial basis of the algebra
$\mathcal{A}_{\mathscr{P}_n}$. The algebra
$\mathcal{A}_{\mathscr{P}_n}$ is are graded. For a graded algebra
$\mathcal{A}=\mathcal{A}^0\oplus \mathcal{A}^1\oplus
\mathcal{A}^2\cdots$, the Hilbert series of the algebra
$\mathcal{A}$, is the formal power series in $x$ given by
$$Hilb\mathcal{A}=\sum\limits_{i\geq 0}x^idim\mathcal{A}^i;$$
there exists a polynomial $P_{\mathcal{A}}(x)$ with rational
coefficients (called the {\it Hilbert polynomial} of $\mathcal{A}$)
such that $P_{\mathcal{A}}(k)=dim\mathcal{A}^i$ for all sufficiently
large $i$.

\begin{thm}\label{HilbertApolynomial} The Hilbert polynomial $P_{\mathcal{A}_{\mathscr{P}_n}}(x)$ of the
algebra $\mathcal{A}_{\mathscr{P}_n}$ is
$x^{\lfloor\frac{1}{2}(n-1)\rfloor}\mathscr{P}_n(\frac{1}{x})$.
\end{thm}
\begin{proof} Note that the number of the multichains
$S_{n,j_1}\preceq S_{n,j_2}\preceq \cdots \preceq S_{n,j_i}$in
$\mathscr{P}_n$ is equal to the dimension of $\mathcal{A}_n^i$.
Corollary \ref{zetefunction} implies that
$\dim\mathcal{A}_n^i=i^{\lfloor\frac{n-1}{2}\rfloor}\mathscr{P}_n(\frac{1}{i})$
for any $i\geq 1$. Since
$deg(\mathscr{P}_n(x))=\lfloor\frac{n-1}{2}\rfloor$, we have
$P_{\mathcal{A}_{\mathscr{P}_n}}(x)=x^{\lfloor\frac{1}{2}(n-1)\rfloor}\mathscr{P}_n(\frac{1}{x})$.
\end{proof}

\begin{thm}\label{HilbertseriesA}For any $n\geq 3$, the Hilbert series $Hilb\mathcal{A}_{\mathscr{P}_n}(x)$ of the algebra
$\mathcal{A}_{\mathscr{P}_n}$ satisfies the following recurrence:
when $n$ is even,
$$Hilb\mathcal{A}_{\mathscr{P}_{n+1}}(x)=xHilb\mathcal{A}'_{\mathscr{P}_{n}}(x)+Hilb\mathcal{A}_{\mathscr{P}_{n}}(x);$$
when $n$ is odd,
\begin{eqnarray*}Hilb\mathcal{A}_{\mathscr{P}_{n+3}}(x)&=&xHilb\mathcal{A}_{\mathscr{P}_{n+2}}'(x)+Hilb\mathcal{A}_{\mathscr{P}_{n+2}}(x)+\frac{4n}{n+3}xHilb\mathcal{A}_{\mathscr{P}_{n+1}}'(x)\\
&&-\frac{8n}{n+3}xHilb\mathcal{A}_{\mathscr{P}_{n}}'(x)-\frac{4n}{n+3}x^2Hilb\mathcal{A}_{\mathscr{P}_{n}}''(x),\end{eqnarray*}
 where the notation ``$~'~$" denotes differentiation of functions,
with the initial condition
$Hilb\mathcal{A}_{\mathscr{P}_{3}}(x)=\frac{1}{(1-x)^2}$ and
$Hilb\mathcal{A}_{\mathscr{P}_{4}}(x)=\frac{1+x}{(1-x)^2}$.\end{thm}
\begin{proof} By Theorem \ref{HilbertApolynomial}, we have
$Hilb\mathcal{A}_{\mathscr{P}_{n}}(x)=1+\sum\limits_{i\geq
1}i^{\lfloor\frac{n-1}{2}\rfloor}\mathscr{P}_n(\frac{1}{i})x^i$. The
results follow Theorem \ref{theoremfpolynomial}.
\end{proof}

In general, we may suppose that
$Hilb\mathcal{A}_{\mathscr{P}_n}(x)=\frac{A_{\mathscr{P}_n}(x)}{(1-x)^{\lceil\frac{1}{2}n\rceil}}$,
where $A_{\mathscr{P}_n}(x)$ is a polynomial.
\begin{cor}For any $n\geq 3$,  $A_{\mathscr{P}_n}(x)$ satisfies the following recurrence:
when $n$ is even,
$$A_{\mathscr{P}_{n+1}}(x)=x(1-x)A_{\mathscr{P}_n}'(x)+\left[(\frac{1}{2}n-1)x+1\right]A_{\mathscr{P}_n}(x);$$
when $n$ is odd,
\begin{eqnarray*}(1-x)A_{\mathscr{P}_{n+3}}(x)&=&x(1-x)A_{\mathscr{P}_{n+2}}'(x)+\frac{(n+1)x+2}{2}A_{\mathscr{P}_{n+2}}(x)+\frac{4n}{n+3}x(1-x)^2A_{\mathscr{P}_{n+1}}'(x)\\
&&+\frac{2n(n+1)}{n+3}x(1-x)A_{\mathscr{P}_{n+1}}(x)-\frac{4n}{n+3}x(1-x)[(n-1)x+2]A_{\mathscr{P}_{n}}'(x)\\
&&-\frac{n(n+1)}{n+3}x[(n-1)x+4]A_{\mathscr{P}_{n}}(x)-\frac{4n}{n+3}x^2(1-x)^2A_{\mathscr{P}_{n}}''(x)\end{eqnarray*}
with the initial conditions $A_{\mathscr{P}_{3}}(x)=1$ and
$A_{\mathscr{P}_{4}}(x)=1+x$.\end{cor}

\begin{proof} By Theorem \ref{HilbertseriesA}, we immediately obtain
the desired results after simple computations. \end{proof}

Note that a sequence $S_{n,i_1}, S_{n,i_2},\cdots, S_{n,i_j}$ of
elements in the poset $\mathscr{P}_n$ forms a chain if and only if
the monomial $x_{i_1}x_{i_2}\cdots x_{i_j}$ is nonvanishing in the
algebra $\mathcal{B}_{\mathscr{P}_n}$.

\begin{thm}For any $n\geq 3$, the Hilbert series $Hilb\mathcal{B}_{\mathscr{P}_n}(x)$ of the algebra
$\mathcal{B}_{\mathscr{P}_n}$ is
$$Hilb\mathcal{B}_{\mathscr{P}_n}(x)=1+\sum\limits_{i=1}^{\lfloor\frac{1}{2}(n-1)\rfloor+1}\sum{n\choose{d_1,d_2,\cdots,d_{i+1}}}\frac{2d_{i+1}-n}{n}x^i,$$
where the second sum is over all $(d_1,\cdots,d_{i+1})$ such that
$\sum\limits_{k=1}^{i+1}d_k=n$, $d_1\geq 0$, $d_k\geq 1$ for all
$2\leq k\leq i$ and $d_{i+1}\geq
n-\lfloor\frac{n-1}{2}\rfloor$.\end{thm} \begin{proof} Note that the
number of the chains $S_{n,j_1}\prec S_{n,j_2}\prec \cdots \prec
S_{n,j_i}$in $\mathscr{P}_n$ is equal to the dimension of
$\mathcal{B}_{\mathscr{P}_n}^i$. Theorem \ref{thmchain} implies that
$\dim\mathcal{B}_{\mathscr{P}_n}^i=\sum{n\choose{d_1,d_2,\cdots,d_{i+1}}}\frac{2d_{i+1}-n}{n}$
for any $i\geq 1$, where the sum is over all $(d_1,\cdots,d_{k+1})$
such that $\sum\limits_{k=1}^{i+1}d_k=n$, $d_1\geq 0$, $d_k\geq 1$
for all $2\leq k\leq i$ and $d_{i+1}\geq
n-\lfloor\frac{n-1}{2}\rfloor$. This completes the proof.\end{proof}

%%%%%%%%%%%%%%%%%%%%%%%

%%%%%%%%%%%%%%%%%%%%%%%%%%%%%%%%%%%%%%%%%%%%%%%%%%%%

\end{document}